\documentclass[reqno,12pt]{amsart}
\usepackage{amsmath,amsthm,amssymb,amsfonts,amscd}
\usepackage{xypic}
\theoremstyle{plain} 
	\newtheorem{thm}{Theorem}[section]
	\newtheorem*{thm*}{Theorem}
	
	\newtheorem{lem}[thm]{Lemma}
	
	\newtheorem{prop}[thm]{Proposition}
	\newtheorem{conj}[thm]{Conjecture}
	\newtheorem*{conj*}{Conjecture}
	
\theoremstyle{definition}
	\newtheorem{defn}[thm]{Definition}

\theoremstyle{remark}
	\newtheorem{rem}[thm]{Remark}
	
	\newtheorem*{pf}{Proof}
\numberwithin{equation}{section}
\def\CC{{\mathbb C}}

\def\LL{{\mathbb L}}

\def\RR{{\mathbb R}}
\def\ZZ{{\mathbb Z}}

\def\D{{\mathcal D}}
\def\E{{\mathcal E}}

\def\L{{\mathcal L}}

\def\N{{\mathcal N}}
\def\O{{\mathcal O}}
\def\P{{\mathcal P}}

\def\T{{\mathcal T}}

\def \mf#1#2#3#4{
\xymatrix{
{#1}\  \ar@<0.4ex>[r]^{{#2}} & \ {#4}
\ar@<0.4ex>[l]^{{#3}}
}
}

\def \mfs#1#2#3#4{\!
\xymatrix@C=1,5em{{#1} \! \ar@<0.2ex>[r]^{{#2}} & \! {#4}
\ar@<0.2ex>[l]^{{#3}}
}
\!}

\def \mfl#1#2#3#4{
\xymatrix@C=2.6em{{#1}\  \ar@<0.4ex>[r]^{{#2}} &\  {#4}
\ar@<0.2ex>[l]^{{#3}}
}
}

\def \mfss#1#2#3#4{\!
\xymatrix@C=1.5em{{#1} \ar@<0.3ex>[r]^{{#2}} & {#4}
\ar@<0.3ex>[l]^{{#3}}
}
\!}
\begin{document}
\title{On entropy for autoequivalences of the derived category of curves}
\date{\today}
\author{Kohei Kikuta}
\address{Department of Mathematics, Graduate School of Science, Osaka University, 
Toyonaka Osaka, 560-0043, Japan}
\email{k-kikuta@cr.math.sci.osaka-u.ac.jp}
\begin{abstract}
To an exact endofunctor of a triangulated category with a split-generator,
the notion of entropy is given by Dimitrov--Haiden--Katzarkov--Kontsevich, 
which is a (possibly negative infinite) real-valued function of a real variable. 
It is important to evaluate the value of the entropy at zero
in relation to the topological entropy.
In this paper, we study the entropy at zero of an autoequivalence of the derived category of a complex smooth projective curve,
and prove that it coincides with the natural logarithm of the spectral radius of the induced automorphism on its numerical Grothendieck group.
\end{abstract}
\maketitle
\section{Introduction}
To an exact endofunctor of a triangulated category with a split-generator, 
the notion of entropy is given by Dimitrov--Haiden--Katzarkov--Kontsevich in \cite{DHKK}. 
It is a (possibly negative infinite) real-valued function of a real variable 
motivated by an analogy with the topological entropy. 
It is known that the topological entropy of a surjective holomorphic endomorphism of a compact K\"{a}hler manifold coincides with the natural logarithm of the spectral radius of the induced action on the cohomology, 
which is the fundamental theorem of Gromov-Yomdin \cite{Gro1,Gro2,Yom} (see also Theorem 3.6 in \cite{Ogu}). 
Concerning the entropy at zero of an exact endofunctor of a complex smooth proper variety, 
there is a lower bound, under a certain technical assumption, by the natural logarithm of the spectral radius of the induced action on the Hochschild homology (Theorem~2.9 in \cite{DHKK}).
If the variety is projective, then the categorical and the topological entropies of 
surjective endomorphisms satisfying the assumption coincide (Theorem 2.12 in \cite{DHKK}). 
It is indeed possible to show this without the technical assumption used in \cite{DHKK}, 
which is given in \cite{KT}.

In this paper, we study the entropy of autoequivalences of the derived category of a smooth projective curve, and prove the following theorem, which is a natural generalization of the fundamental theorem of Gromov-Yomdin.
\newpage
\begin{thm}[Theorem~\ref{main-theorem}]
Let $C$ be a complex smooth projective curve and 
$F$ an autoequivalence of the bounded derived category $\D^b(C)$ of coherent sheaves on $C$.  
The entropy $h(F)$ coincides with the natural logarithm of the spectral radius $\rho([F])$ of the induced automorphism $[F]$ on the numerical Grothendieck group $\N(C)$ of $\D^b(C)$. 
In particular, $\rho([F])$ is an algebraic number.
\end{thm}

The contents of this paper is as follows. 
In Section 2, we recall the definition and some basic properties of the entropy of exact endofunctors by \cite{DHKK}.
We shall study the entropy for an autoequivalence of the bounded derived category of coherent sheaves on a smooth projective curve in Section 3.
The non-trivial case is when a curve is an elliptic curve. 
We shall introduce the notion of the autoequivalence of type-{\bf m} (Definition \ref{type-m}), which behaves very well in two points:
it gives a representative of a conjugacy class of the automorphism group of the numerical Grothendieck group preserving the Euler form (Proposition \ref{Z-conjugate}), 
and it enables us to compute explicitly the entropy (Proposition \ref{computable}). 
In Section 4, we shall recall the notion of LLS-period and give a proof of
Proposition \ref{Z-conjugate},
the key proposition for our main theorem.

\begin{sloppypar}
{\bf Acknowledgements}.\  
I am grateful to my supervisor, Professor Atsushi Takahashi for guidances, supports and encouragements. 
He suggested to us Conjecture~\ref{conj} motivated by the theorem of Gromov--Yomdin and Theorems~2.9, 2.12 in \cite{DHKK}.  
I would like to thank Yuuki Shiraishi for giving many comments on this paper. 

\end{sloppypar}
\section{Preliminaries}
\subsection{Notations and terminologies}
Throughout this paper, we work over the base field $\CC$ and all
triangulated categories are $\CC$-linear
and not equivalent to the zero category. The translation functor on a
triangulated category is denoted by $[1]$.

A triangulated category $\T$ is called {\it split-closed} if every
idempotent in $\T$ splits, namely, if it contains all direct summands of
its objects, and it is called {\it thick} if it is split-closed and closed
under isomorphisms.
For an object $M\in\T$, we denote $\langle M\rangle$ by the smallest thick
triangulated subcategory containing $M$.
An object $G\in\T$ is called a {\it split-generator} if $\langle
G\rangle=\T$.
A triangulated category $\T$ is said to be {\it of finite type} if for all
$M,N\in\T$ we have
$\sum_{n\in\ZZ}\dim_\CC{\rm Hom}_{\T}(M,N[n])<\infty$.
\newpage
\subsection{Complexity}
From now on, $\T$ denotes a triangulated category of finite type.
\begin{defn}[Definition~2.1 in \cite{DHKK}]
For each $M,N\in \T$, define the function $\delta_{\T,t}(M,N):\RR\longrightarrow\RR_{\geq0}\cup\{ \infty \}$ in $t$ by
{\small 
\begin{equation*}
\delta_{\T,t}(M,N):= 
\begin{cases}
0
 & \text{ if }N\cong0\\
\inf\left\{
\displaystyle\sum_{i=1}^p {\rm exp}(n_i t)~
\middle
|~
\begin{xy}
(0,5) *{0}="0", (20,5)*{A_{1}}="1", (30,5)*{\dots}, (40,5)*{A_{p-1}}="k-1", (60,5)*{N\oplus N'}="k",
(10,-5)*{M[n_{1}]}="n1", (30,-5)*{\dots}, (50,-5)*{M[n_{p}]}="nk",
\ar "0"; "1"
\ar "1"; "n1"
\ar@{.>} "n1";"0"
\ar "k-1"; "k" 
\ar "k"; "nk"
\ar@{.>} "nk";"k-1"
\end{xy}
\right\}
 & \text{ if }N\in\langle M\rangle \\
\infty
 & \text{ if }N\not\in\langle M\rangle.
\end{cases}
\end{equation*}
}
The function $\delta_{\T,t}(M,N)$ is called the {\it complexity} of $N$ with  respect to $M$.
\end{defn}
\begin{rem}
If $\T$ has a split-generator $G$ and $M\in\T$ is not isomorphic to a zero object, 
then an inequality $1\leq\delta_{\T,0}(G,M)<\infty$ holds.
\end{rem} 
We recall some basic properties of the complexity.  
\begin{lem}[Proposition~2.3 in \cite{DHKK}]\label{complexity}
Let $M_1,M_2,M_3\in\T$.
\begin{enumerate}
\item
If $M_1\cong M_3$, then $\delta_{\T,t}(M_1,M_2)=\delta_{\T,t}(M_3,M_2)$. 
\item
If $M_2\cong M_3$, then $\delta_{\T,t}(M_1,M_2)=\delta_{\T,t}(M_1,M_3)$. 
\item
If $M_2\not\cong 0$, then $\delta_{\T,t}(M_1,M_3)\leq\delta_{\T,t}(M_1,M_2)\delta_{\T,t}(M_2,M_3)$.
\item
Let $\T'$ be a triangulated category of finite type.
We have $\delta_{\T',t}(F(M_1),F(M_2))\leq\delta_{\T,t}(M_1,M_2)$ for any exact functor $F:\T\longrightarrow \T'$.
\end{enumerate}
\end{lem}
\begin{lem}[Lemma~2.4 in \cite{DHKK}]\label{per(k)}
Let $\D^b(\CC)$ be the bounded derived category of finite dimensional $\CC$-vector spaces. 
For $M\in\D^b(\CC)$, we have the following inequality 
\begin{equation}\label{eq:2}
\delta_{\D^b(\CC),t}(\CC,M)=\sum_{l\in\ZZ} \left(\dim_{\CC}H^{l}(M)\right)\cdot e^{-lt}.
\end{equation}
\end{lem}
\subsection{Entropy of endofunctors}
\begin{defn}[Definition~2.5 in \cite{DHKK}]
Let $G$ be a split-generator of $\T$ and $F$ an exact endofunctor of $\T$ such that $F^n G\not\cong 0$ for $n\ge 0$. 
The {\it entropy} of $F$ is the function $h_t(F):\RR\longrightarrow\{ -\infty \}\cup\RR$ given by
\begin{equation*}
h_{t}(F):=\displaystyle\lim_{n\rightarrow\infty}\frac{1}{n}\log \delta_{\T,t}(G,F^{n}G).
\end{equation*}
\end{defn}
It follows from Lemma~2.6 in \cite{DHKK} that the entropy is well-defined.
\begin{lem}\label{two-split-gen}
Let $G,G'$ be split-generators of $\T$ and $F$ an exact endofunctor of $\T$ such that $F^n G,F^nG'\not\cong 0$ for $n\ge 0$. 
Then
\[
h_t(F)=\displaystyle\lim_{n\rightarrow\infty}\frac{1}{n}\log \delta_{\T,t}(G,F^{n}G').
\]
\end{lem}
\begin{pf}
It follows from Lemma~\ref{complexity} (iii) and (iv) that
\[
\delta_{\T,t}(G,F^nG')\leq \delta_{\T,t}(G,F^nG)\delta_{\T,t}(F^nG,F^nG')\leq \delta_{\T,t}(G,G')\delta_{\T,t}(G,F^nG).
\]
Similarly, we have $\delta_{\T,t}(G,F^nG)\leq\delta_{\T,t}(G',G)\delta_{\T,t}(G,F^nG')$, which yields the statement.
\qed
\end{pf}
\begin{lem}[Section~2 in \cite{DHKK}]\label{entropy}
Let $G$ be a split-generator of $\T$ and $F_1,F_2$ exact endofunctors of $\T$ such that $F_i^n G\not\cong 0$ for $i=1,2$ and $n\ge 0$. 
\begin{enumerate}
\item
If $F_{1} \cong F_{2}$, then $h_t(F_{1})=h_t(F_{2})$. 
\item
We have $h_t(F_{1}^m)= mh_t(F_{1})$ for $m\geq1$.
\item
We have $h_t([m])=mt$ for $m\in\ZZ$. 
\item
If $F_1F_2\cong F_2F_1$, then $h_t(F_1F_2)\leq h_t(F_1)+h_t(F_2)$.
\item
If $F_1=F_2[m]$ for $m\in\ZZ$, then $h_t(F_1)=h_t(F_2)+mt$. 
\end{enumerate}
\end{lem}
\begin{lem}
Let $G$ be a split-generator of $\T$ and $F_1,F_2$ exact endofunctors of $\T$ such that $F_i^n G\not\cong 0$ for $i=1,2$ and $n\ge 0$. 
Then we have $h_t(F_1F_2)=h_t(F_2F_1)$.
\end{lem}
\begin{pf}
It follows from Lemma~\ref{complexity} (iii) and (iv) that
\begin{eqnarray*}
\delta_{\T,t}(G,(F_1F_2)^nG)
&=& \delta_{\T,t}(G,F_1(F_2F_1)^{n-1}F_2G)\\
&\leq& \delta_{\T,t}(G,F_1G)\delta_{\T,t}(F_1G,F_1(F_2F_1)^{n-1}F_2G)\\
&\leq& \delta_{\T,t}(G,F_1G)\delta_{\T,t}(G,(F_2F_1)^{n-1}G)\\
&\cdot& \delta_{\T,t}((F_2F_1)^{n-1}G,(F_2F_1)^{n-1}F_2G)\\
&\leq& \delta_{\T,t}(G,F_1G)\delta_{\T,t}(G,F_2G)\delta_{\T,t}(G,(F_2F_1)^{n-1}G),
\end{eqnarray*}
which gives $h_t(F_1F_2)\le h_t(F_2F_1)$.
Similarly, we can show $h_t(F_2F_1)\le h_t(F_1F_2)$, which yields the statement. 
\qed
\end{pf}
\begin{lem}
Let $F$ be an exact endofunctor of $\T=\langle G\rangle$ such that $F^nG\not\cong 0$ for $n\ge 0$, $F'$ an autoequivalence of $\T$, and $F'^{-1}$ the quasi-inverse of $F'$. 
The entropy is a class function, namely, 
$h_t(F'FF'^{-1})=h_t(F)$. 
\end{lem}
\begin{pf}
The inequality $h_t(F'FF'^{-1})\leq h_t(F)$ follows from 
\begin{eqnarray*}
\delta_{\T,t}(G,(F'FF'^{-1})^nG)  
&=& \delta_{\T,t}(G,F'F^nF'^{-1}G)\\
&\leq& \delta_{\T,t}(G,F'F^nG)\delta_{\T,t}(F'F^nG,F'F^nF'^{-1}G)\\
&\leq& \delta_{\T,t}(G,F'G)\delta_{\T,t}(F'G,F'F^nG)\delta_{\T,t}(G,F'^{-1}G)\\
&\leq& \{ \delta_{\T,t}(G,F'G)\delta_{\T,t}(G,F'^{-1}G) \} \delta_{\T,t}(G,F^nG)
\end{eqnarray*}
Considering $F'':=F'^{-1}FF'$, we have the reverse inequality. 
\qed
\end{pf}
\subsection{The cases of smooth projective varieties}
Let $X$ be a smooth projective variety.
The bounded derived category $\D^b(X)$ of coherent sheaves over $X$ is of finite type 
and it has a split-generator. In particular, $G:=\oplus_{i=1}^{\dim X+1}\O_X(i)$ and its dual $G^*:=\oplus_{i=1}^{\dim X+1}\O_X(-i)$ 
are known to be split-generators (cf. Theorem 4 in \cite{Orl}). 
The following proposition enables us to compute entropies. 
\begin{prop}[cf. Theorem~2.7 in \cite{DHKK}]\label{smooth-proper}
Let $G,G'$ be split-generators of $\D^b(X)$ and $F$ an autoequivalence of $\D^b(X)$. 
The entropy $h_t(F)$ is given by
\begin{equation*}
h_t(F)=\lim_{n\rightarrow\infty}\frac{1}{n}\log\delta'_{\D^b(X),t}(G,F^nG'),
\end{equation*}
where 
\begin{equation*}
\delta'_{\D^b(X),t}(M,N):=\sum_{m\in\ZZ}\left(\dim_{\CC} {\rm Hom}_{\D^b(X)}(M,N[m])\right)\cdot e^{-mt},\quad M,N\in\D^b(X).
\end{equation*}
\end{prop}
\begin{pf}
The following is proven in the proof of Theorem~2.7 in \cite{DHKK}.
\begin{lem}\label{delta-delta'}
There exist $C_1(t),C_2(t)$ for $t\in\RR$ such that
\[
C_1(t)\delta_{\D^b(X),t}(G,M)
\leq\delta'_{\D^b(X),t}(G,M)
\leq C_2(t)\delta_{\D^b(X),t}(G,M),\quad M\in\D^b(X).
\]
In particular, for each $M\in \D^b(X)$ we have
\begin{equation}
\lim_{n\rightarrow\infty}\frac{1}{n}\log\delta_{\D^b(X),t}(G,M)
=\lim_{n\rightarrow\infty}\frac{1}{n}\log\delta'_{\D^b(X),t}(G,M).
\end{equation}
\end{lem}
Together with Lemma~\ref{two-split-gen}, we have
\[
h_t(F)=\lim_{n\rightarrow\infty}\frac{1}{n}\log\delta_{\D^b(X),t}(G,F^nG')
=\lim_{n\rightarrow\infty}\frac{1}{n}\log\delta'_{\D^b(X),t}(G,F^nG').
\]
We finished the proof of the proposition.
\qed
\end{pf}
In order to study the structure of the entropy we prepare some terminologies.
For $M,N\in\D^b(X)$, set 
\begin{equation}
\chi(M,N):=\sum_{n\in\ZZ}(-1)^n\dim_\CC{\rm Hom}_{\D^b(X)}(M,N[n]).
\end{equation}
It naturally induces a bilinear form on the Grothendieck group $K_0(X)$ of $\D^b(X)$, called the {\it Euler form}, 
which is denoted by the same letter $\chi$. 
Then {\it numerical Grothendieck group} $\N(X)$ is defined as the quotient of $K_0(X)$ by the radical of $\chi$ 
(which is well-defined by the Serre duality). 
It is known that $\N(X)$ is a free abelian group of finite rank by Hirzebruch-Riemann-Roch theorem. 

The group of isomorphism classes of autoequivalences of $\D^b(X)$ is denoted by 
${\rm Auteq}(\D^b(X))$. 
We denote the group of isomorphisms of $\N(X)$ preserving $\chi$ by ${\rm Aut}_\ZZ(\N(X), \chi)$. 
Then we have the group homomorphism
\[
{\rm Auteq}(\D^b(X))\rightarrow{\rm Aut}_\ZZ(\N(X), \chi),\quad F\mapsto [F].
\]
\begin{defn}
For an $F\in{\rm Auteq}(\D^b(X))$, the {\it spectral radius} $\rho([F])$ is 
the maximum of absolute values of eigenvalues of the induced endomorphism $[F]\in {\rm Aut}_\ZZ(\N(X), \chi)$.
\end{defn}
From now on, we shall only consider entropies at $t=0$. 
For simplicity, set $\delta_{\D^b(X)}:=\delta_{\D^b(X),0}$, $h(F):=h_0(F)$ and so on. 
\begin{defn}
For an $f\in {\rm Aut}(X)$, the {\it categorical entropy} $h_{cat}(f)$ is the entropy $h(f^*)$ of its derived functor at $t=0$. 
\end{defn}
It is known that the topological entropy $h_{top}(f)$ of a surjective holomorphic endomorphism $f\in {\rm End}(M)$ of a compact K\"{a}hler manifold $M$ coincides with the natural logarithm of the spectral radius $\rho'(f)$ of the induced action on the cohomology, 
which is the fundamental theorem of Gromov-Yomdin \cite{Gro1,Gro2,Yom} (see also Theorem 3.6 in \cite{Ogu}). 
Concerning the entropy $h(F)$ of an exact endofunctor $F$ of $\D^b(X)$ of a complex smooth proper variety $X$, 
there is a lower bound, under a certain technical assumption, by the natural logarithm of the spectral radius of the induced action on the Hochschild homology (Theorem~2.9 in \cite{DHKK}). 
If $X$ is projective, $h_{cat}(f)=\log\rho'(f)$ for a surjective endomorphism $f\in {\rm End}(X)$ satisfying the assumption (Theorem 2.12 in \cite{DHKK}). 
It is indeed possible to show that $h_{cat}(f)=\log\rho([f^*])=\log\rho'(f)=h_{top}(f)$ without the technical assumption used in \cite{DHKK}, which is given in \cite{KT}.

From the above, it is natural to expect the following conjecture in general. 
\begin{conj}\label{conj}
Let $X$ be a smooth projective variety.
For each $F\in {\rm Auteq}(\D^b(X))$, 
the entropy should coincide with the natural logarithm of the spectral radius:
\begin{equation}
h(F)=\log\rho([F]).
\end{equation}
\end{conj}
The conjecture is true in the case of the ample canonical or anti-canonical sheaf (\cite{KT}). 
\section{Entropy for curves}
The main result of this paper is the following. 
\begin{thm}\label{main-theorem}
Let $C$ be a smooth projective curve.
For each $F\in{\rm Auteq}(\D^b(C))$, the entropy coincides with 
the natural logarithm of the spectral radius:
\begin{equation}
h(F)=\log\rho([F]).
\end{equation}
\end{thm}
\subsection{Standard autoequivalences}
\begin{defn}
The {\it standard autoequivalence group} is the subgroup of ${\rm Auteq}(\D^b(X))$ given by
\begin{equation}
{\rm Auteq}^{st}(\D^b(X)):={\rm Aut}(X)\ltimes({\rm Pic}(X)\times\ZZ[1]).
\end{equation}
\end{defn}
\begin{prop}\label{standard-curve}
Let $C$ be a smooth projective curve. 
For each standard autoequivalence $F\in{\rm Auteq}^{st}(\D^b(C))$, we have $h(F)=0$. 
\end{prop}
\begin{pf}
Each standard autoequivalence $F$ is represented as $F=f^*(-\otimes\L)[m]$. 
Since Lemma~\ref{entropy} (v) gives $h(F)=h(f^*(-\otimes\L))$, we can assume $F=f^*(-\otimes\L)$. 
Let $\O_C(1)$ be a very ample invertible sheaf and set $G:=\O_C(1)\oplus\O_C(2)$.
Recall that $G$ and $G^*=\O_C(-1)\oplus\O_C(-2)$ are split-generators of $\D^b(C)$. 
By Lemma~\ref{delta-delta'}, we have 
\[
h(F)=\lim_{n\rightarrow\infty}\frac{1}{n}\log\delta_{\D^b(C)}(G,F^n G^*)
=\lim_{n\rightarrow\infty}\frac{1}{n}\log\delta'_{\D^b(C)}(G,F^nG^*).
\]
If $\deg\L>0$, then $G^*\otimes F^nG^*$ becomes a sum of invertible sheaves with positive degrees for large $n$, 
and hence we have ${\rm Hom}_{\D^b(X)}(G,F^nG^*[i])=0$ for $i\ne 0$ and large $n$.
If $\deg\L\le 0$, then $G^*\otimes F^nG^*$ becomes a sum of invertible sheaves with negative degrees for $n\ge 0$, 
and hence we have ${\rm Hom}_{\D^b(X)}(G,F^nG^*[i])=0$ for $i\ne 1$ and $n\ge 0$.
Since $f^*$ acts trivially on $\N(C)$, we obtain
\[
\lim_{n\rightarrow\infty}\frac{1}{n}\log\delta'_{\D^b(C)}(G,F^nG^*)
=\lim_{n\rightarrow\infty}\frac{1}{n}\log\left|\chi(G,F^nG^*)\right|
=\lim_{n\rightarrow\infty}\frac{1}{n}\log\left|\chi(G,G^*\otimes \L^n)\right|.
\]
It follows from the polynomial-growth of $\chi(G,G^*\otimes\L^n)$ (cf. Lemma 2.14 in \cite{DHKK}) that
\[
h(F)=\lim_{n\rightarrow\infty}\frac{1}{n}\log\left|\chi(G,G^*\otimes\L^n)\right|=0.
\]
We finished the proof of the proposition.
\qed
\end{pf}
Let $C$ be a smooth projective curve which is not an elliptic curve. 
Then it follows from \cite{BO} that ${\rm Auteq}(\D^b(C))={\rm Auteq}^{st}(\D^b(C))$.
Therefore, for the proof of our main theorem we only need to consider the case of an elliptic curve. 
\subsection{Elliptic curves}
Let $X=(E,x_0)$ be an elliptic curve with a closed point $x_0\in E$. 
The numerical Grothendieck group $\N(E)$ has a canonical basis given by $[\O_E]$ and $[\O_{x_0}]$ and 
the Euler form $\chi$ with respect to this basis is given by 
\[
\begin{pmatrix}
0&1\\
-1&0
\end{pmatrix},
\]
which yields the isomorphism ${\rm Aut}_\ZZ(\N(E), \chi)\cong {\rm {\rm SL}}(2,\ZZ)$. 
Consider the Fourier-Mukai functor $\Phi_\P:=\RR p_{2*}(\LL p_1^*(-)\otimes^\LL\P)$, 
an autoequivalence of $\D^b(E)$ first introduced by Mukai \cite{Muk}, 
where $\P\in {\rm coh}(E\times E)$ is the Poincar\'{e} bundle and $p_1$ (resp. $p_2$) is the projection 
$E\times E\longrightarrow E$ to the first (resp. second) component.
Set $S:=\Phi_\P,T:=-\otimes\O_E(x_0)$. 
These satisfy the following relations \cite{Muk,ST}:
\begin{equation}\label{relations-with-S,T}
S^2\cong(-1)^*[-1],~(TS)^3\cong S^2.
\end{equation}
With a canonical basis of $\N(E)$ given by $[\O_E]$ and $[\O_{x_0}]$, we have 
\begin{equation}\label{ST}
[S]=
\begin{pmatrix}
0&1\\
-1&0
\end{pmatrix},
\quad [T]=
\begin{pmatrix}
1&0\\
1&1
\end{pmatrix}.
\end{equation}
It is a well-known fact that the natural map 
$\phi: {\rm Auteq}(\D^b(E))\longrightarrow {\rm Aut}_\ZZ(\N(E), \chi)$ is surjective and 
gives the following exact sequence
\begin{equation}
\{1\}\longrightarrow {\rm Aut}(E)\ltimes({\rm Pic}^0(E)\times\ZZ[2])\longrightarrow 
{\rm Auteq}(\D^b(E))\longrightarrow {\rm {\rm SL}}(2,\ZZ)\longrightarrow \{1\}.
\end{equation}
\begin{lem}\label{ker-phi}
The map $h:{\rm Auteq}(\D^b(E))\longrightarrow \RR_{\ge 0}$, $F\mapsto h(F)$ factors through ${\rm {\rm SL}}(2,\ZZ)$.
\end{lem}
\begin{pf}
Set $\O_E(1):=\O_E(3x_0)$ and $G:=\O_E(1)\oplus\O_E(2)$. 
Suppose that an element $F\in {\rm Auteq}(\D^b(E))$ is of the form $F=F'F_1$ with 
$F'\in {\rm Auteq}(\D^b(E))$ and $F_1\in {\rm Aut}(E)\ltimes({\rm Pic}^0(E)\times\ZZ[2])$. 
Then there exist $F_2,\dots, F_n\in {\rm Aut}(E)\ltimes({\rm Pic}^0(E)\times\ZZ[2])$ such that 
$F^n=(F'F_1)^n=F'^nF_n \cdots  F_1$. We have 
\[
\delta_{\D^b(E)}(G,F^n G^*)
\le \delta_{\D^b(E)}(G,F'^n G)\delta_{\D^b(E)}(G,F_n\cdots F_1 G^*),
\]
and hence, 
\[
h(F)\le h(F')+\lim_{n\rightarrow\infty}\frac{1}{n}\log\delta_{\D^b(E)}(G,F_n\cdots F_1 G^*).
\]
Since $F_i$ is of the form $f_i^*(-\otimes\L_i)[2l_i]$ for some $f_i\in {\rm Aut}(E)$, $\L_i\in {\rm Pic}^0(E)$ and $l_i\in\ZZ$, $G^*\otimes F_n\cdots F_1 G^*[-2l_n-\cdots -2l_1]$ is a sum of anti-ample invertible sheaves.
Therefore, we have $\delta'_{\D^b(E)}(G,F_n\cdots F_1 G^*)=\left|\chi(G,F_n\cdots F_1 G^*)\right|$.
Since $F_i$ acts trivially on $\N(E)$, we have $\left|\chi(G,F_n\cdots F_1 G^*)\right|=\left|\chi(G,G^*)\right|$, 
which implies 
\begin{multline*}
\lim_{n\rightarrow\infty}\frac{1}{n}\log\delta_{\D^b(E)}(G,F_n\cdots F_1 G^*)
=\lim_{n\rightarrow\infty}\frac{1}{n}\log\delta'_{\D^b(E)}(G,F_n\cdots F_1 G^*)\quad (\text{by Lemma~\ref{delta-delta'}})\\
=\lim_{n\rightarrow\infty}\frac{1}{n}\log\left|\chi(G,F_n\cdots F_1 G^*)\right|
=\lim_{n\rightarrow\infty}\frac{1}{n}\log\left|\chi(G,G^*)\right|=0.
\end{multline*}
Therefore $h(F)\leq h(F')$. 
We also have $h(F')\leq h(F)$ since $F'=FF_1^{-1}$ and $F_1^{-1}$ belongs to  ${\rm Aut}(E)\ltimes({\rm Pic}^0(E)\times\ZZ[2])$. 
\qed
\end{pf}
\begin{lem}
The entropies of $S,T,TS$ are all equal to zero.
\end{lem}
\begin{pf}
From the relations~(\ref{relations-with-S,T}) and Lemma~\ref{entropy}~(ii) $h(S)=0$ and $h(TS)=0$ follow, 
$h(T)=0$ is given by Proposition~\ref{standard-curve}.
\qed
\end{pf}
The characteristic polynomial of a matrix $A\in {\rm SL}(2,\ZZ)$ is given by
\begin{equation}
x^2-({\rm tr}A)x+1=0. 
\end{equation}
If $|{\rm tr}[F]|\leq2$, then $\rho(A)=1$. 
If ${\rm tr}[F]>2$, then eigenvalues of $A$ are $\rho,\rho^{-1}$, where $\rho$ is a real number greater than 1. 
In particular, the spectral radius $\rho(A)$ of $A$ is $\rho$, which is an algebraic number.
Since $h(F)=h(S^2F)$ and $[S^2F]=-[F]$ for any $F\in{\rm Auteq}(\D^b(E))$, 
we do not need to consider the case of ${\rm tr}[F]<-2$. 
\begin{prop}\label{leq2}
For every autoequivalence $F\in{\rm Auteq}(\D^b(E))$ with $|{\rm tr}[F]|\leq2$, 
\[
h(F)=0=\log\rho([F]).
\] 
\end{prop}
\begin{pf}
First we consider the case of ${\rm tr}[F]=-1,0,1$. 
From the form of characteristic polynomial of $[F]$, we have
\[
[F]^2=-
\left(
\begin{matrix}
1&0\\
0&1
\end{matrix}
\right)
\text{ or }
[F]^3=\pm
\left(
\begin{matrix}
1&0\\
0&1
\end{matrix}
\right).
\]
From Lemma~\ref{entropy} (ii), clearly their entropies are equal to 0. 
Then in the case of $|{\rm tr}[F]|=2$, 
representatives of conjugacy classes are given as follows:
\[
\left(
\begin{matrix}
1&0\\
0&1
\end{matrix}
\right),~
\left(
\begin{matrix}
1&n\\
0&1
\end{matrix}
\right)
=[(TST)^n]~(n>0).
\]
We have
\begin{eqnarray*}
h((TST)^n)
&=& h((S^{-1}T^{-1}S^{-1})^n)\quad(\text{by (\ref{relations-with-S,T})})\\
&=& h(S^{-1}T^{-1}S^{-2}T^{-1}\cdots T^{-1}S^{-2}T^{-1}S^{-1})\\
&=& h(S^{-2}T^{-1}S^{-2}T^{-1}\cdots T^{-1}S^{-2}T^{-1})\\
&=& h(T^{-n})=nh(T^{-1})=0. 
\end{eqnarray*}
Clearly $\rho([F])=1$ in both cases. 
We finished the proof of the proposition.  
\qed
\end{pf}
We shall introduce the notion of the autoequivalence of type-{\bf m}, which behaves very well in two points:
it gives a representative of a conjugacy class of ${\rm SL}(2,\ZZ)$, 
and it enables us to compute explicitly the entropy. 
\begin{defn}\label{type-m}
Let ${\bf m}=(m_{2n},\dots,m_1)$, $n\ge 1$ be an ordered sequence of positive integers. 
An autoequivalence $F\in{\rm Auteq}(\D^b(E))$ is {\it of type-{\bf m}} if $F$ is isomorphic to following types:
\[
F_{\bf m}=
\begin{cases}
S^2T^{m_{2n}}ST^{-m_{2n-1}}S\cdots T^{m_2}ST^{-m_1}S
 & \text{if $n$:odd}\\
T^{m_{2n}}ST^{-m_{2n-1}}S\cdots T^{m_2}ST^{-m_1}S
 & \text{if $n$:even}.
\end{cases}
\]
\end{defn}
\begin{prop}\label{Z-conjugate}
For any $F\in{\rm Auteq}(\D^b(E))$ with ${\rm tr}[F]>2$, 
there exists an autoequivalence $F_{\bf m}$ of type-{\bf m} such that $[F]$ is conjugate to $[F_{\bf m}]$. 
\end{prop}
Since the proof is elementary but technical, we shall give it in Appendix. 
\begin{prop}\label{computable}
For each $F\in{\rm Auteq}(\D^b(E))$ of type-{\bf m}, we have $h(F)=\log\rho([F])$. 
\end{prop}
\begin{pf}
First, we show the following lemma. 
\begin{lem}\label{basic-behavior}
Let $F$ be an autoequivalence of type-{\bf m}.
For each indecomposable locally free sheaf $\E$ on $E$ with positive degree, 
$F(\E)$ is an indecomposable locally free sheaf with positive degree up to even translations. 
\end{lem}
\begin{pf}
Note that each autoequivalence sends an indecomposable object to an indecomposable one.
It follows from equation~\eqref{ST} that if $\E$ is an indecomposable locally free sheaf with non-zero degree then
$S(\E)$ is a locally free sheaf up to translations since the abelian category ${\rm coh}(E)$ is hereditary. 
Obviously, $T$ sends a locally free sheaf to a locally free sheaf.
In addition, for positive numbers $r_1,d_1$, we have
\[
\begin{pmatrix}
r_2\\
d_2
\end{pmatrix}:=
[T^{-m_{2i-1}}S]
\begin{pmatrix}
r_1\\
d_1
\end{pmatrix}=
\begin{pmatrix}
0&1\\
-1&-m_{2i-1}
\end{pmatrix}\begin{pmatrix}
r_1\\
d_1
\end{pmatrix}=
\begin{pmatrix}
d_1\\
-r_1-m_{2i-1}d_1
\end{pmatrix}
\in
\begin{pmatrix}
\ZZ_{>0}\\
\ZZ_{<0}
\end{pmatrix},
\]
\[
\begin{pmatrix}
r_3\\
d_3
\end{pmatrix}:=
[T^{m_{2i}}S]
\begin{pmatrix}
r_2\\
d_2
\end{pmatrix}=
\begin{pmatrix}
0&1\\
-1&m_{2i}
\end{pmatrix}\begin{pmatrix}
r_2\\
d_2
\end{pmatrix}=
\begin{pmatrix}
d_2\\
-r_2+m_2d_2
\end{pmatrix}\in
\begin{pmatrix}
\ZZ_{<0}\\
\ZZ_{<0}
\end{pmatrix},
\]
\[
\begin{pmatrix}
r_4\\
d_4
\end{pmatrix}:=
[T^{-m_{2i+1}}S]
\begin{pmatrix}
r_3\\
d_3
\end{pmatrix}
=
\begin{pmatrix}
0&1\\
-1&-m_{2i+1}
\end{pmatrix}
\begin{pmatrix}
r_3\\
d_3
\end{pmatrix}
=
\begin{pmatrix}
d_3\\
-r_3-m_{2i+1}d_3
\end{pmatrix}\in
\begin{pmatrix}
\ZZ_{<0}\\
\ZZ_{>0}
\end{pmatrix},
\]
\[
\begin{pmatrix}
r_5\\
d_5
\end{pmatrix}:=
[T^{m_{2i+2}}S]
\begin{pmatrix}
r_4\\
d_4
\end{pmatrix}
=
\begin{pmatrix}
0&1\\
-1&m_{2i+2}
\end{pmatrix}
\begin{pmatrix}
r_4\\
d_4
\end{pmatrix}
=
\begin{pmatrix}
d_4\\
-r_4+m_{2i+2}d_4
\end{pmatrix}\in
\begin{pmatrix}
\ZZ_{>0}\\
\ZZ_{>0}
\end{pmatrix},
\]
which implies that $F(\E)$ is a locally free sheaf with positive degree up to even translations. 
\qed
\end{pf}
For any positive integer $l$, $F^l\O_E(1),F^l\O_E(2)$ are indecomposable locally free sheaves with positive degrees up to even translations by Lemma~\ref{basic-behavior}. 
Hence, it easily follows from Atiyah's results (cf. Lemma~1.1 in \cite{Har}) that 
${\rm Hom}_{\D^b(E)}(G^*,F^lG[n])=0$ for $l\ge 0$ and any odd integer $n$, which gives
\begin{eqnarray*}
h(F)
&=& \lim_{l\rightarrow\infty}\frac{1}{l}\log\delta_{\D^b(E)}(G^*,F^l G) =\lim_{l\rightarrow\infty}\frac{1}{l}\log\delta'_{\D^b(E)}(G^*,F^l G)\quad (\text{by Lemma~\ref{delta-delta'}})\\
&=& \lim_{l\rightarrow\infty}\frac{1}{l}\log\chi(G\otimes F^l G).
\end{eqnarray*}
We obtain $\lim_{l\rightarrow\infty}\sqrt[l]{\chi(G\otimes F^l G)}=\rho([F])$, 
in this case by some elementary computations of linear algebra, which finished the proof of Proposition~\ref{computable}
\qed
\end{pf}
By Propositions \ref{leq2}, \ref{Z-conjugate} and \ref{computable}, we obtain the following
\begin{thm}
For each $F\in{\rm Auteq}(\D^b(E))$, we have $h(F)=\log\rho([F])$. 
\end{thm}
To summarize, we finished the proof of our main theorem, Theorem~\ref{main-theorem}. 
\section{Appendix}
\subsection{LLS-period}
This subsection is based on Chapter~7 in \cite{Kar}. 
\begin{defn}[Definition~7.12 in \cite{Kar}]
Let $A$ be a matrix in ${\rm SL}(2,\ZZ)$ represented as
\[
A=
\begin{pmatrix}
a&b\\
c&d
\end{pmatrix}. 
\] 
If $0< a\leq c<d$ , we call $A$ a {\it reduced matrix}.  
\end{defn}
\begin{rem}
The definition of a reduced matrix is slightly different from the one in \cite{Kar}, 
where $A$ is a reduced matrix if and only if $0\leq a\leq c<d$. 
A reduced matrix in the sense of \cite{Kar} with $a=0$ is represented as 
\[
\begin{pmatrix}
0&-1\\
1&d
\end{pmatrix}
=
\begin{pmatrix}
1&-1\\
0&1
\end{pmatrix}
\begin{pmatrix}
1&d-2\\
1&d-1
\end{pmatrix}
\begin{pmatrix}
1&1\\
0&1
\end{pmatrix}
,\quad d>2, 
\]
namely, it is conjugate to a reduce matrix in the sense of this paper.
Thus, for the proof of Proposition~\ref{Z-conjugate}, it is sufficient to consider only the case of $a>0$. 
\end{rem}
\begin{prop}[Theorem~7.13 in \cite{Kar}]
For any matrix $A\in {\rm SL}(2,\ZZ)$ with ${\rm tr}A>2$, 
there exists a reduced matrix $A'$ such that $A$ is conjugate to $A'$. 
\end{prop}
For each $A\in {\rm SL}(2,\ZZ)$ with ${\rm tr}A>2$, 
one can associate a cyclically ordered sequence $LLS(A)$ of an even length with positive integers, called the {\it LLS-period}.  
For the details, see Definition~7.10 in \cite{Kar}. 
\begin{prop}[Theorem~7.14 in \cite{Kar}]
Let $A\in {\rm SL}(2,\ZZ)$ be a reduced matrix. 
Suppose that there exists positive integers $a_1,\dots, a_{2n-1},a_{2n}$, $n\ge 1$ such that 
\[
\frac{c}{a}=[a_1;a_2,\dots, a_{2n-1}]:=
a_1+\frac{1}{a_2+\frac{1}{\ddots +\frac{1}{a_{2n-2}+\frac{1}{a_{2n-1}}}}},\quad 
\lfloor\frac{d-1}{c}\rfloor=a_{2n},
\]
where $[a_1;a_2,\dots, a_{2n-1}]$ is the canonical continued fraction representation of $c/a$ and 
$\lfloor\frac{d-1}{c}\rfloor$ is the integer part of $(d-1)/c$.
Then the LLS-period $LLS(A)$ of $A$ is given by
\begin{equation}
LLS(A)=(a_1,a_2,\cdots,a_{2n-1},a_{2n}). 
\end{equation}
\end{prop}
\begin{prop}[Proposition~7.11 in \cite{Kar}]\label{complete-invariant}
For reduced matrices $A,A'\in {\rm SL}(2,\ZZ)$ with ${\rm tr}A,{\rm tr}A'>2$, 
$A$ is conjugate to $A'$ if and only if $LLS(A)$ coincides with $LLS(A')$. 
\end{prop}
\subsection{The proof of Proposition~\ref{Z-conjugate}}
Recall that ${\bf m}=(m_{2n},\dots,m_1)$, $n>1$, where $m_i$ are all positive integers. 
\begin{prop}[Proposition~\ref{Z-conjugate}]
For each $F\in{\rm Auteq}(\D^b(E))$ with ${\rm tr}[F]>2$, 
there exists $F_{\bf m}\in{\rm Auteq}(\D^b(E))$ of type-{\bf m} such that $[F]$ is conjugate to $[F_{\bf m}]$. 
\end{prop}
\begin{pf}
The following lemma is essential. 
\begin{lem}\label{LLS-of-typem}
For each $F_{\bf m}\in{\rm Auteq}(\D^b(E))$ of type-{\bf m}, 
$[F_{\bf m}]$ is a reduced matrix such that $LLS([F_{\bf m}])={\bf m}$.
\end{lem}
\begin{pf}
Set
\[
[F_{\bf m}]=
\begin{pmatrix}
a_n&b_n\\
c_n&d_n
\end{pmatrix},~
\alpha_n:=b_n-m_1a_n,~
\beta_n:=d_n-m_1c_n. 
\]
We shall show the following five properties by induction on $n$: 
(i) $0<a_n$, (ii) $a_n\leq c_n$, (iii) $\displaystyle\frac{c_n}{a_n}=[m_{2n};m_{2n-1},\dots,m_2]$, 
(iv) $0\leq\alpha_n<a_n$, (v) $0<\beta_n\leq c_n$. 

If $n=1$, then all the properties follows from the computation 
\[
[S^2T^{m_2}ST^{-m_1}S]
=
\begin{pmatrix}
a_1&b_1\\
c_1&d_1
\end{pmatrix}
=
\begin{pmatrix}
1&m_1\\
m_2&m_1m_2+1
\end{pmatrix}. 
\]

Suppose that the properties holds for $n$ and we shall show them for $n+1$. 
The property (i) is clear by
\begin{eqnarray*}
\begin{pmatrix}
a_{n+1}&b_{n+1}\\
c_{n+1}&d_{n+1}
\end{pmatrix}
&=&
\begin{pmatrix}
1&m_{2n+1}\\
m_{2n+2}&m_{2n+1}m_{2n+2}+1
\end{pmatrix}
\begin{pmatrix}
a_n&b_n\\
c_n&d_n
\end{pmatrix}\\
&=& 
\begin{pmatrix}
a_n+c_nm_{2n+1}&b_n+d_nm_{2n+1}\\
a_nm_{2n+2}+c_n(m_{2n+1}m_{2n+2}+1)&b_nm_{2n+2}+d_n(m_{2n+1}m_{2n+2}+1)
\end{pmatrix}. 
\end{eqnarray*}
The properties (ii) and (iii) follow from the computation 
\begin{equation*}
\frac{c_{n+1}}{a_{n+1}}
=\frac{a_nm_{2n+2}+c_n(m_{2n+1}m_{2n+2}+1)}{a_n+c_nm_{2n+1}}
=m_{2n+2}+\frac{c_n}{a_n+c_nm_{2n+1}}
=m_{2n+2}+\frac{1}{m_{2n+1}+\frac{1}{\frac{c_n}{a_n}}}. 
\end{equation*}
The property (iv) holds since we have 
\begin{eqnarray*}
b_{n+1}
&=& b_n+d_nm_{2n+1}\\
&=& (m_1a_n+\alpha_n)+d_nm_{2n+1}\\
&=& (m_1a_n+m_1m_{2n+1}c_n)-m_1m_{2n+1}c_n+\alpha_n+d_nm_{2n+1}\\
&=& m_1a_{n+1}+(\alpha_n+m_{2n+1}(d_n-m_1c_n))\\
&=& m_1a_{n+1}+(\alpha_n+m_{2n+1}\beta_n)\\
&<& m_1a_{n+1}+(a_n+m_{2n+1}c_n)\\
&=& m_1a_{n+1}+a_{n+1}. 
\end{eqnarray*}
Finally, we see the property (v) by
\begin{eqnarray*}
d_{n+1}
&=& b_nm_{2n+2}+d_n(m_{2n+1}m_{2n+2}+1)\\
&=& (m_1a_n+\alpha_n)m_{2n+2}+(m_1c_n+\beta_n)(m_{2n+1}m_{2n+2}+1)\\
&=& m_1c_{n+1}+\alpha_nm_{2n+2}+\beta_n(m_{2n+1}m_{2n+2}+1)\\
&\leq& m_1c_{n+1}+\{a_nm_{2n+2}+c_n(m_{2n+1}m_{2n+2}+1)\}\\
&=& m_1c_{n+1}+c_{n+1}. 
\end{eqnarray*}
Consequently, all the properties hold for $n+1$ and hence they do for all positive integers. 

From the properties (i),(ii),(v) and the definition of $\beta_n$, we have $0<a_n\leq c_n<d_n$, which means that $[F_{\bf m}]$ is a reduced matrix. 
By the equality
\begin{eqnarray*}
\frac{d_n-1}{c_n}
&=& m_1+\frac{\beta_n-1}{c_n}
\end{eqnarray*}
and the property (v), we have
\[
LLS([F_{\bf m}])=(m_{2n},\dots,m_1)={\bf m}. 
\]
We finished the proof of the lemma.
\qed
\end{pf}
Lemma~\ref{LLS-of-typem} and Proposition~\ref{complete-invariant} yield the statement of the proposition. 
\qed
\end{pf}

\end{document}